\documentclass[12pt]{article}
 \usepackage[german]{babel}
\usepackage{amssymb}
\usepackage{amsmath}
\usepackage{graphicx}
\usepackage{bbm}
\usepackage[right]{eurosym}
\usepackage[left=2cm,right=2cm,top=2cm,bottom=3cm]{geometry}
\def\vs{\vspace}
\def\noi{\noindent}

\def\IN{\mathbb N}

\def\IR{\mathbb R}
\def\IC{\mathbb C}

\def\an{\mathrm{an}}
\def\exp{\mathrm{exp}}

\def\ma{\mathcal}
\begin{document}
\begin{center}
{\bf \Large Piecewise Weierstra\ss $ $ preparation and division for}
\end{center}

\begin{center}
{\bfseries\Large o-minimal holomorphic functions}
\end{center}

\vs{0.5cm}
\centerline{Tobias Kaiser}

\vspace{0.7cm}\noi \footnotesize {{\bf Abstract.} 
Given an o-minimal structure expanding the field of reals, we show a piecewise Weierstrass preparation theorem and a piecewise Weierstrass  division theorem for definable holomorphic functions. In the semialgebraic setting and for the structure of globally subanalytic sets and functions we obtain the corresponding results for definable real analytic functions.
As an application we show a definable global Nullstellensatz for principal ideals.}

\normalsize
\section*{Introduction}

The Weierstrass preparation theorem and the Weierstrass division theorem are the key tool for local complex analytic geometry (see for example Gunning and Rossi [11, Chapter II \& III]and \L ojasiewicz [15, Chapters I-III]). They are also used
for example in the proof of the important Oka's coherence theorem (see [11, Chapter IV]).\\
We deal with Weierstrass theorems in the o-minimal setting. Classes of functions allowing  Weierstrass preparation and division have been one of the main sources to establish and to analyze o-minimal structures (see Denef and Van den Dries [4], Van den Dries [5], Van den Dries and Speissegger [8, 9], and [14]). They also play a role in the recent applications of o-minimality to number theory (see Binyamini and Novikov [1,2]).\\
In [19], Peterzil and Starchenko have given for o-minimal structures on the reals strong versions of classical results of complex analytic geometry as the Remmert-Stein theorem and Remmert's proper mapping theorem. In a series of papers [16, 17, 18] they have developed complex analysis in an arbitrary o-minimal expansion of a real closed field. In [17, Section 2.4] they establish in this setting the classical, i.e. local, Weierstrass  theorems. In [18, Section 11] they prove a finite version of the coherence theorem where they use in a special setting globally prepared functions. \\
We show in general for an arbitrary o-minimal expansion of the reals\footnote{The results hold in o-minimal expansions of arbitrary real closed fields. We prefer the presentation in the real setting to keep it at a reasonable length.}  a finite Weierstrass preparation and a finite Weierstrass division theorem for definable holomorphic functions.
Finite means that the domain of definition can be covered by finitely many definable open sets such that we have there, after some coordinate transformation, a global  preparation resp. division theorem. Our setting includes for example the case of complex Nash functions (see for example Fatabbi and Tancredi [10] for the latter).\\
The strategy for establishing the results is the following: The first observation is that, given a nonzero definable holomorphic function on a domain, the order of its power series expansions at the points of the domain is uniformely bounded.
After applying a linear coordinate transformation from a suitable finite list one can assume that these power series expansions are piecewise regular of constant order with respect to the last variable. 

\rule{16cm}{0.01cm}

{\footnotesize{\itshape 2010 Mathematics Subject Classification:} 03C64, 32B05, 32B20, 32C07}
\newline
{\footnotesize{\itshape Keywords and phrases:} O-minimal structure, Weierstra\ss $ $ preparation and division, Nullstellensatz}

Along such a piece, the data from the classical Weierstrass preparation resp. division are definable in the parameters. Using o-minimal geometric arguments, we can glue these definable families to obtain piecewise Weierstrass preparation and division.\\
Applying the result from [12] that real analytic functions definable in the structure $\IR_\an$ have a definable holomorphic extension, we get in this o-minimal structure piecewise Weierstrass preparation and division for definable real analytic functions. The same holds for semialgebraic real analytic functions; i.e. for Nash functions.\\
As an application we establish a definable global Nullstellensatz for principal ideals (compare with [15, Chapter II, \S 5]).

\section*{Notations}

By $|X|$ we denote the number of elements of a finite set $X$.\\
We identify $\IC^n$ with $\IR^{2n}$ in the canonical way.
Given a set $Z\subset \IC^p\times \IC^n$ and a function $f:Z\to \IC$ we write
$Z_a:=\big\{z\in \IC^n\colon (a,z)\in Z\big\}$ and $f_a:Z_a\to \IC, z\mapsto f(a,z)$ for $a\in\IC^p$.\\
By $|z|$ we denote the euclidean norm and by $||z||$ the sum norm of $z\in \IC^n$.
We write $z=(z_1,\ldots,z_n)=(z^\dagger,z_n)$. For $X\subset \IC^n$ we set $X^\dagger:=\pi(X)$ where $\pi:\IC^n\to \IC^{n-1}$ denotes the projection onto the first $n-1$ components.\\
For $a\in \IC^n$ and $r\in \IR_{>0}^n$ let
$\Delta(a,r)=\Delta^n(a,r):=\big\{z\in\IC^n\colon |z_j-a_j|<r_j\mbox{ for all }j\big\}$ be the open multidisc with center $a$ and multiradius $r$.\\
A function $u:X\to\IC$ is called a {\bf unit} if $u(x)\neq 0$ for all $x\in X$. The {\bf degree} of a polynomial is abbreviated by $\mathrm{deg}$. A polynomial is called {\bf monic} if its leading coefficient equals $1$.\\
Given a holomorphic function $f:U\to \IC$ on an open subset $U$ of $\IC^n$ and $a\in U$ we denote by $\mathrm{ord}_a(f)\in \IN_0\cup\{\infty\}$ the {\bf order} of the power series expansion of $f$ at $a$.
The function $f$ is {\bf regular of order $k$ in the variable $z_n$ at $a$} if
the power series expansion of $f(a_1,\ldots,a_{n-1},z_n)$ at $a_n$ has order $k$.\\
A domain in $\IC^n$ is a nonempty and connected open subset of $\IC^n$.

\vs{0.5cm}
\noindent
Below, we fix an o-minimal structure $\ma{M}$ expanding the field of reals. ``Definable'' means ``definable in $\ma{M}$ with parameters'' unless otherwise stated.

\section{Preparations}

{\bf 1.1 Definition}
	\begin{itemize}
		\item[(a)]
		Let $U$ be an open subset of $\IC^n$ that is definable.
		By $\ma{O}^{\mathrm{def}}(U)=\ma{O}_\ma{M}^{\mathrm{def}}(U)$ we denote the ring of functions $f:U\to \IC$ that are holomorphic and definable.
		\item[(b)] Let $U$ be an open subset of $\IR^n$ that is definable.
		By $\ma{C}^{\omega,\mathrm{def}}(U)=\ma{C}_\ma{M}^{\omega,\mathrm{def}}(U)$ we denote the ring of functions $f:U\to \IR$ that are real analytic and definable.
	\end{itemize}

\vs{0.2cm}
\noindent
The following definition and remark can be stated for arbitrary holomorphic functions. We formulate it in the definable case.  

\newpage
{\bf 1.2 Definition}

\vs{0.1cm}
	Let $U$ be an open and definable set in $\IC^n$ and let $f\in\ma{O}^\mathrm{def}(U)$.
	We say that $f$ is polynomial in the variable $z_n$ of degree $k$ if $f\in \ma{O}^\mathrm{def}(U^\dagger)[z_n]$ and $\mathrm{deg}(f)=k$, i.e. if there are $a_0,\ldots,a_k\in \ma{O}^\mathrm{def}(U^\dagger)$ with $a_k\neq 0$ such that
	$$f(z)=a_0(z^\dagger)+a_1(z^\dagger)z_n+\ldots +a_k(z^\dagger)z_n^k.$$

\vs{0.5cm}

{\bf 1.3 Proposition}

\vs{0.1cm}
	Let $U$ be an open and definable set in $\IC^n$ and let $f\in \ma{O}^\mathrm{def}(U)$. Assume that $U_{z^\dagger}$ is connected for all $z^\dagger\in U^\dagger$. Then the following are equivalent:
	\begin{itemize}
		\item[(i)] $f$ is polynomial in the variable $z_n$ of degree at most $k$,
		\item[(ii)] $\partial^{k+1}f/\partial z_n^{k+1}$ vanishes on $U$. 
	\end{itemize}
{\bf Proof:}

\vs{0.1cm}
	The direction (i) $\Rightarrow$ (ii) is clear.
	For the direction (ii) $\Rightarrow$ (i) let $g:=\partial^k f/\partial z_n^k$. Then
	$\partial g/\partial z_n$ vanishes on $U$. Given $z^\dagger\in U^\dagger$ we obtain that the function $g_{z^\dagger}:U_{z^\dagger}\to\IC$ is constant since $U_{z^\dagger}$ is connected by assumption. Hence we see that $g\in \ma{O}^\mathrm{def}(U^\dagger)$. 
	Let $h:=f-g z_n^k/k!\in\ma{O}^\mathrm{def}(U)$. Then $\partial^k h/\partial z_n^k$ vanishes on $U$. By induction we get the claim.
\hfill$\Box$

\vs{0.5cm}
\noindent 
The above holds similarly in the real analytic setting.

\vs{0.2cm}
\noindent The following result which we will use later has been mentioned in Van den Dries and Miller [7, Remark on p. 1376]. For the readers' convenience we formulate the proof.

\vs{0.5cm}
{\bf 1.4 Remark}

\vs{0.1cm}
	Let $U$ be a definable and open set in $\IR^n$. Then $C^{\omega,\mathrm{def}}(U)$ is weakly noetherian in the sense of Tougeron [20], i.e. the following conditions hold:
	\begin{itemize}
		\item[(1)] $C^{\omega,\mathrm{def}}(U)$ contains the real polynomials in $n$ variables,
		\item[(2)] $C^{\omega,\mathrm{def}}(U)$ is stable under taking derivatives, 
		\item[(3)] every descending sequence of zero sets of functions in $C^{\omega,\mathrm{def}}(U)$ is ultimately stationary. 
	\end{itemize}
{\bf Proof:}

\vs{0.1cm}
	The ring $C^{\omega,\mathrm{def}}(U)$ is stable under taking derivatives and clearly contains the polynomials. 
	It fulfils condition b) on [20, p. 824]:
	For example, take the function $x\mapsto 1/(1+|x|^2)$.
	Moreover, by o-minimality, it fulfills condition c'') on [20, p. 831]:
	A finite dimensional linear subspace of $C^{\omega,\mathrm{def}}(U)$ is a definable family of functions. We obtain the claim by the uniform finiteness property (see for example Van den Dries [6, Chapter 3, Corollary (3.7)]). 
	Hence by [20, Th\'{e}or\`{e}me II'] the ring $C^{\omega,\mathrm{def}}(U)$ is topologically noetherian which implies that the ring is weakly noetherian (see [20, p. 824]). 
\hfill$\Box$

\vs{0.5cm} 
The following observation can be proved in several ways (see [19, Theorem 7.1] and the successive discussion there). For the readers' convenience we include a proof.

\vs{0.5cm}
{\bf 1.5 Proposition}

\vs{0.1cm}	
	Let $\Omega$ be a definable domain in $\IC^n$ and let $f\in\ma{O}^{\mathrm{def}}(\Omega)$. Assume that $f\neq 0$.	
	Then there is some $K\in \IN_0$ such that $\mathrm{ord}_a(f)\leq K$ for every $a\in \Omega$.
	
\vs{0.1cm}
{\bf Proof}

\vs{0.1cm}	
	For $k\in \IN_0$ let
	$$\omega_k:\Omega\to \IR_{\geq 0},a\mapsto \big\vert \big(D^\alpha f(a)\big)_{||\alpha||\leq k}\big\vert^2.$$
	Then $\omega_k\in C^{\omega,\mathrm{def}}(\Omega)$ for every $k\in \IN_0$ (where $\Omega$ is viewed as a subset of $\IR^{2n}$).
	For $k\in \IN_0$ let $\ma{Z}(\omega_k):=\big\{a\in \Omega\colon \omega_k(a)=0\big\}$ be the zero set of $\omega_k$. The sequence $\big(\ma{Z}(\omega_k)\big)_{k\in\IN_0}$ is clearly decreasing.
	By Remark 1.4 this sequence is ultimately stationary.
	Hence there is some $K\in \IN_0$ such that
	$\bigcap_{k\in\IN_0}\ma{Z}(\omega_k)=
	\bigcap_{k=0}^K\ma{Z}(\omega_k)$.
	Since $\Omega$ is a domain and $f$ is a nonzero holomorphic function on $\Omega$ we see that the former set is empty.
	This gives the claim.
\hfill$\Box$

\vs{0.5cm}
{\bf 1.6 Definition}

\vs{0.1cm}	
	In the above situation we call
	$$\mathrm{ord}(f):=\max\big\{\mathrm{ord}_a(f)\colon a\in \Omega\big\}$$
	the {\bf order of} $f$.

\section{Classical Weierstrass theorems with definable parameters}

\noindent 
We establish in the definable setting parameterized versions of Weierstrass preparation and division. They are obtained by a common model theoretic argument from the fact that by 
[17, Section 2.4] the classical Weierstrass theorems hold in o-minimal expansions of arbitrary real closed fields.

\vs{0.2cm}
\noindent
Let $A\subset \IC^p$ be definable.

\vs{0.5cm}
{\bf 2.1 Definition}
	\begin{itemize}
		\item[(a)]
		By $\mathfrak{U}^n(A)$ we denote the set of definable subsets $\ma{X}$ of $\IC^p\times \IC^n$ such that $\ma{X}_a$ is an open neighbourhood of $0$ in $\IC^n$ for every $a\in A$.
		\item[(b)] 
		A function $F:\ma{X}\to \IC$ is called {\bf dph}
		(definably parameterized holomorphic) if $\ma{X}\in \mathfrak{U}^n(A)$ and $F_a:\ma{X}_a\to A$ is holomorphic for every $a\in A$.
	\end{itemize}

\vs{0.2cm}
{\bf 2.2 Remark}
\begin{itemize}
		\item[(1)]
		Let $\varphi:A\to \IR_{>0}^n$ be definable. Then
		$$\ma{X}(A,\varphi):=\bigcup_{a\in A}\Big(\{a\}\times \Delta\big(0,\varphi(a)\big)\Big)\in \mathfrak{U}^n(A).$$
		\item[(2)]
		Let $\ma{X}\in \mathfrak{U}^n(A)$. Then there is a definable function $\varphi:A\to\IR_{>0}^n$ such that $\ma{X}(A,\varphi)\subset \ma{X}$.
	\end{itemize}

\vs{0.2cm}
{\bf 2.3 Proposition}

\vs{0.1cm}	
	Let $\ma{X}\in \mathfrak{U}^n(A)$ and let $F:\ma{X}\to \IC$ be a dph function. Assume that there is some $k\in \IN_0$ such that, for every $a\in A$, the holomorphic function $F_a$ is regular of order $k$ in the variable $z_n$ at $0$.
	Then we find $\ma{Y}\in \mathfrak{U}^n(A)$ with $\ma{Y}\subset \ma{X}$ and dph functions
	$H,U:Y\to \IC$ such that the following holds for every $a\in A$:
	\begin{itemize}
		\item[(1)] $F_a=H_aU_a$ on $\ma{Y}_a$,
		\item[(2)] $H_a$ is a Weierstrass polynomial of degree $k$ in the variable $z_n$ at $0$,
		\item[(3)] $U_a$ is a unit on $\ma{Y}_a$.
	\end{itemize}
{\bf Proof:}

\vs{0.1cm}
	Assume that the statement does not hold.
	Let $\Sigma$ be the set of all triples $(\ma{Y},H,U)$ where $\ma{Y}\in \mathfrak{U}^n(A)$ with $\ma{Y}\subset \ma{X}$ and $H,U:\ma{Y}\to \IC$ are dph-functions such that $H_a$ is a Weierstrass polynomial of degree $k$ in the variable $z_n$ at $0$ and $U_a$ is a unit on $\ma{Y}_a$ for all $a\in A$.
	Let $\ma{L}$ be the language which has a relation symbol for every definable subset in any dimension and view $\ma{M}$ as an $\ma{L}$-structure in the natural way.
	For $\sigma=(\ma{Y},H,U)\in \Sigma$ let
	$\varphi_\sigma(a)$ be a $\ma{L}$-formula expressing that $a\in A$ and $F_a=H_aU_a$.
	By assumption we get that the set of $\ma{L}$-formulas 
	$\Xi(a):=\big\{\neg \varphi_\sigma(a)\colon \sigma\in \Sigma\big\}$
	in the variable $a$ is finitely realizable.
	Hence there is an elementary extension $\ma{N}$ of $\ma{M}$ such that the type $\Xi(a)$ is realized in $\ma{N}$.
	Let $R$ be the universe of $\ma{N}$ and let $K:=R[\sqrt{-1}]$. One identifies $K$ with $R^2$ in the canonical way.
	So there is some $b\in K^p$ such that $\ma{N}\models\neg\varphi_\sigma(b)$ for all $\sigma\in\Sigma$.
	We have for the canonical liftings of $A$, $\ma{X}$ and $F$ to $\ma{N}$
	that $b\in A$, $\ma{X}_b$ is an open subset of $K^n$ containing $0$, $F_b$ is $K$-differentiable in the sense of [17] and $F_b$ is regular of order $k$ in $z_n$ at $0$. By 
	[17, Theorem 2.20] we find an open subset $Y_b$ of $K^n$ containing $0$ and $K$-differentiable functions $h_b,u_b$ such that $F_b=h_bu_b$ on $Y_b$, $h_b$ is a Weierstrass polynomial of degree $k$ in the variable $z_n$ at $0$ and $u_b$ is a unit on $Y_b$. 
	There is some $\sigma=(\ma{Y},H,U)\in\Sigma$ such that $\ma{Y}_b=Y_b,H_b=h_b$ and $U_b=u_b$. So $\ma{N}\models \varphi_\sigma(b)$, contradiction.
\hfill$\Box$

\vs{0.5cm}
In the same way we obtain the following.

\vs{0.5cm}
{\bf 2.4 Proposition}

\vs{0.1cm}
	Let $\ma{X}\in \mathfrak{U}^n(A)$ and let
	$F:\ma{X}\to \IC$ be a dph function. Assume that there is some $k\in \IN_0$ such that, for every $a\in A$, the function $F_a$ is regular of order $k$ in the variable $z_n$ at $0$.
	Let $G:\ma{X}\to \IC$ be a dph function. Then there are $\ma{Y}\in \mathfrak{U}^n(A)$ with $\ma{Y}\subset \ma{X}$ and dph functions $Q:\ma{Y}\to \IC, R:\ma{Y}\to \IC$ such that
	$$G=QF+R$$ and $R_a$ is polynomial in the variable $z_n$ of degree at most $k-1$ for all $a\in A$.

\newpage
\section{Piecewise definable Weierstrass theorems}

Recall, that by o-minimality, a definable set has only finitely many connected components and that each of them is definable. So we work in the following with domains.

\vs{0.2cm} 
\noindent
Let $\Omega$ be a definable domain in $\IC^n$ and let $f\in\ma{O}^{\mathrm{def}}(\Omega)$ with $f\neq 0$.
For $a\in \Omega$ and $k\in \IN_0$ let
$$J_a^k:\IC^n\to \IC, w\mapsto\sum_{||\alpha||=k}D^\alpha f(a)w^\alpha.$$
We have that
$$f(z)=\sum_{k=0}^\infty J_a^k(z-a)$$
on open multidiscs centered at $a$.
For $k\in \big\{0,\ldots,\mathrm{ord}(f)\big\}$ let
$$\Omega_k:=\big\{a\in  \Omega\colon \mathrm{ord}_a(f)=k\big\}.$$ Note that $\Omega_0$ is a dense open subset of $\Omega$ and that $f|_{\Omega_0}$ is a unit.
We find $\ma{X}\in \mathfrak{U}^n(\Omega)$ such that $a+\ma{X}_a\subset \Omega$ for all $a\in \Omega$.

\vs{0.5cm}
{\bf 3.1 Lemma}

\vs{0.1cm}	
	Let $A\subset \IR^r\times\IR^s$ be definable such that
	$\dim A_x<s$ for all $x\in \IR^r$. Then there is $q\in \IN$ and there are $y_1,\ldots,y_q\in \IR^s\setminus\{0\}$ such that for every
	$x\in \IR^r$ there is some $p\in\{1,\ldots,q\}$ with $y_p\notin A_x$.
	
\vs{0.1cm}
{\bf Proof:}

\vs{0.1cm}
	Choose a decomposition $\ma{C}$ of $A$ into definable cells (see for example [6, Chapter 3]).
	Let $C\in \ma{C}$ and let $\iota_C\in \{0,1\}^{r+s}$ be such that $C$ is a
	$\iota_C$-cell.
	By assumption there is
	$t_C\in\{1,\ldots,s\}$ such that the $t_C^{\mathrm{th}}$-component of $\iota_C$ equals
	$0$.
	For $t\in \{1,\ldots,s\}$ let $\ma{C}_t:=\{C\in \ma{C}\colon t_C=t\}$ and $A_t:=\bigcup_{C\in\ma{C}_t}C$. Assuming for notational reasons that $t=s$ for a moment we have that $(A_t)_{(x,y^\dagger)}$ contains at most $|\ma{C}_t|$ points for every $(x,y^\dagger)\in \IR^r\times\IR^{s-1}$.
	Let $\alpha_t:=|\ma{C}_t|+1$. Then any $\alpha_t$-tuple in $\IR^s\setminus\{0\}$ with pairwise distinct $t$-components do the job for $A_t$ instead of $A$.
	Let $\alpha:=\max\{\alpha_1,\ldots,\alpha_s\}$ and $q:=\alpha |\ma{C}|$. Choose $y_1=(y_{11},\ldots,y_{1s}),\ldots,y_q=(y_{q1},\ldots,y_{qs})\in \IR^s\setminus\{0\}$ such that 
	$y_{p_1t}\neq y_{p_2t}$ for all $p_1,p_2\in \{1,\ldots,q\}$ with $p_1\neq p_2$ and all $t\in\{1,\ldots,s\}$. We show that $y_1,\ldots,y_q$ do the job: Let $x\in \IR^r$. Assume that $y_1,\ldots,y_q\in A_x$. By the definition of $q$ we find $C\in \ma{C}$ and pairwise distinct indices $j_1,\ldots,j_\alpha\in \{1,\ldots,q\}$ such that $y_{j_1},\ldots,y_{j_\alpha}\in C_x$. Let $t\in\{1,\ldots,s\}$ with $C\in\ma{C}_t$. Since $\alpha\geq \alpha_t$ and $y_{j_1},\ldots,y_{j_\alpha}$ have pairwise distinct $t$-components we obtain from the above that $y_{j_l}\notin C_x$ for some $l\in\{1,\ldots,\alpha\}$, contradiction.
\hfill$\Box$

\vs{0.5cm}
For $k\in \big\{0,\ldots,\mathrm{ord}(f)\big\}$ let
$$\Gamma_k:=\Big\{(a,w)\in \Omega_k\times\IC^n\colon J_a^k(w)=0\Big\}.$$

\vs{0.5cm}
{\bf 3.2 Proposition}

\vs{0.1cm}	
	Let $k\in \big\{0,\ldots,\mathrm{ord}(f)\big\}$. There is $q_k\in\IN$ and there are linear coordinate transformations $c_{k,1},\ldots,c_{k,q_k}$ such that for every $a\in \Omega_k$ there is $p\in\{1,\ldots,q_k\}$ such that
	$f\circ c_{k,p}$ is regular of order $k$ in the variable $z_n$ at $c_{k,p}^{-1}(a)$.
	
\vs{0.1cm}
{\bf Proof:}

\vs{0.1cm}
	The set
	$\Gamma_k$
	is definable.
	For every $a\in \Omega_k$ we have that the complex dimension of $\Gamma_{k,a}$ is smaller then $n$ since $J_a^k$ is a nonzero polynomial.
	By Lemma 3.1 we find $q_k\in\IN$ and $y_{k,1},\ldots,y_{k,q_k}\in \IC^n\setminus\{0\}$ such that for every $a\in \Omega_k$ there is $p\in\{1,\ldots,q_k\}$ with $J_a^k(y_{k,p})\neq 0$.
	For $p\in \{1,\ldots,q_k\}$ choose $c_{k,p}\in\mathrm{GL}(n,\IC)$ whose last column consists of $y_{k,p}$. Then $c_{k,1},\ldots,c_{k,q_k}$ do the job by [11, Chapter I, Section B, Lemma 2].
\hfill$\Box$

\vs{0.5cm}
\noindent Let
$F:\ma{X}\to \IC, F(a,w)=f(a+w)$.
Then $F$ is a dph function.
Given $a\in \Omega$ we have that
$$F_a(w)=\sum_{k=0}^\infty J_a^k(w)$$
for all $w$ in open multidiscs centered at $0$.

\vs{0.5cm}
{\bf 3.3 Theorem}

\vs{0.1cm}
	Let $k\in\big\{0,\ldots,\mathrm{ord}(f)\big\}$ and let $A$ be a definable subset of $\Omega_k$ such that $f$ is regular in the variable $z_n$ at every $a\in A$. Then there is a finite covering $\ma{V}$ of $A$ by definable open sets that are contained in $\Omega$ such that for every $V\in\ma{V}$ the following holds: There are
	functions 
	$h_V\in \ma{O}^{\mathrm{def}}(V^\dagger)[z_n]$ of degree $k$ in the variable $z_n$ that is monic and
	$u_V\in \ma{O}^{\mathrm{def}}(V)$ which is a unit such that
	$f=h_Vu_V$ on $V$. Moreover, $h_a$ is a Weierstrass polynomial
	of degree $k$ in the variable $z_n$ at every $a\in A$.

\vs{0.1cm}
{\bf Proof:}

\vs{0.1cm}
	Let $\ma{X}_A:=\ma{X}\cap (A\times\IC^n)\in \mathfrak{U}^n(A)$. The dph function $F|_{\ma{X}_A}$ fulfills the assumptions of Proposition 2.3.
	Hence there are  $\ma{Y}\in \mathfrak{U}^n(A)$ with $\ma{Y}\subset \ma{X}_A$ and dph functions
	$H,U:\ma{Y}\to \IC$ such that the following holds for every $a\in A$:
	\begin{itemize}
		\item[(1)] $F_a=H_aU_a$ on $\ma{Y}_a$,
		\item[(2)] $H_a$ is a Weierstrass polynomial of degree $k$ in the variable $z_n$ at $0$,
		\item[(3)] $U_a$ is a unit on $\ma{Y}_a$.
	\end{itemize}
	By Remark 2.2 we may assume that there is a definable function $\varphi:A\to \IR_{>0}^n$ such that $\ma{Y}=\ma{X}(A,\varphi)$. So $\ma{Y}_a=\Delta\big(0,\varphi(a)\big)$ for all $a\in A$. We write $\Theta(a):=\Delta\big(a,\varphi(a)\big)$.
	For $a\in A$ define the holomorphic functions
	$$h^a:\Theta(a)\to \IC, z\mapsto H_a(z-a),$$
	$$u^a:\Theta(a)\to \IC, z\mapsto U_a(z-a).$$
	Then $f=h^au^a$ on $\Theta(a)$ since for $z\in \Theta(a)$
	$$h^a(z)u^a(z)=H_a(z-a)U_a(z-a)=F(a,z-a)=f(a+(z-a))=f(z).$$
	Moreover, $h^a$ is a Weierstrass polynomial of degree $k$ in the variable $z_n$ at $a$ and $u^a$ is a unit. 
	
	\vs{0.2cm}\noindent
	{\bf Claim 1:}
	Let $a\in A$. Let $b\in A$ with $b\in\Theta(a)$. Then $h^a=h^b$ and $u^a=u^b$ on $\Theta(a)\cap\Theta(b)$. 
	
	\vs{0.1cm}\noindent
	{\bf Proof of Claim 1:}
	Let $\Lambda:=\Theta(a)\cap\Theta(b)$. Then $b\in\Lambda$. Developing $h^a$ in a power series at $b$ we obtain that
	$h^a$ is a polynomial of degree $k$ in the variable $z_n$ at $b$.
	Since $f$ is regular of order $k$ in the variable $z_n$ at $b$ and since
	$u^a$ is a unit we have that $h^a$ is a Weierstrass polynomial at $b$. By the uniqueness property of the classical Weierstrass preparation theorem we get that
	the germ of $h^a$ coincides with the germ of $h^b$ at $b$ and that the germ 
	of $u^a$ coincides with the germ of $u^b$ at $b$. Since $\Lambda$ is connected and contains $b$ we obtain the claim.
	\hfill$\Box_{\mathrm{Claim\,1}}$
	
	\vs{0.2cm}\noindent
	Let $Z:=\bigcup_{a\in A}\Theta(a)$. Then $Z$ is an open and definable neighbourhood of $A$ that is contained in $\Omega$.
	Choose a cell decomposition $\ma{C}$ of $A$ such that $\varphi|_C$ is continuous for all $C\in \ma{C}$. Fix $C\in\ma{C}$.
	Being a cell, $C$ is locally closed. Hence we find definable sets $X_C$ and $B_C$ such that $X_C$ is open, $B_C$ is closed and $C=X_C\cap B_C$.
	Let
	$Z_C:=X_C\cap Z$.
	Then $Z_C$ is an open and definable neighbourhood of $C$ that is contained in $\Omega$. Moreover, $C$ is closed in $Z_C$.
	After shrinking $Z_C$ if necessary we find by Van den Dries [6, Chapter 8, Corollary (3.9)] a continuous and definable retraction $R=R_C:Z_C\to C$.
	Since $R|_C$ is the identity and $\varphi|_C$ is continuous we may moreover assume, after shrinking $Z_C$ once more, that
	$z\in \Theta\big(R(z)\big)$ for all $z\in Z_C$.
	Define 
	$$h_C:Z_C\to \IC, z\mapsto h^{R(z)}(z),$$
	$$u_C:Z_C\to \IC, z\mapsto u^{R(z)}(z).$$ 
	
	\vs{0.2cm}\noindent
	{\bf Claim 2:} Let $w\in Z_C$ and let $b:=R(w)$. Then there is a neighbourhood $D$ of $w$ contained in $Z_C\cap\Theta(b)$ such that $h_C|_D=h^b|_D$ and $u_C|_D=u^b|_D$ for all $z\in D$. 
	
	\vs{0.1cm}\noindent
	{\bf Proof of Claim 2:}
	By the continuity of $R$ and of $\varphi|_C$ we find a neighbourhood $D$ of $w$ contained in $Z_C\cap\Theta(b)$ such that $b\in\Theta\big(R(z)\big)$ for all $z\in D$.
	We are done by Claim 1. 
	\hfill$\Box_{\mathrm{Claim\,2}}$
	
	\vs{0.2cm}\noindent	
	Claim 2 implies immediately that $u_C$ and $h_C$ are holomorphic.
	Since $C$ is a cell we have that $C_{z^\dagger}$ is connected for every $z^\dagger\in\IC^{n-1}$. Shrinking $Z_C$ once more if necessary we may assume that 
	$Z_{C,z^\dagger}$ is connected for every $z^\dagger\in\IC^{n-1}$. 
	
	\vs{0.2cm}\noindent
	{\bf Claim 3:} $h_C$ is polynomial of degree $k$ in the variable $z_n$ and monic.
	
	\vs{0.1cm}\noindent
	{\bf Proof of Claim 3:}
	Let $w\in Z_C$. By Claim 2 there is an open neighbourhood $D$ of $w$ contained in $Z_C\cap \Theta(b)$ such that
	$h_C|_D=h^b|_D$ where $b:=R(w)$.
	We have that 
	$\partial^{k+1} h^b/\partial z_n^{k+1}$ vanishes on $D$ since $h^b$ is a Weierstrass polynomial of degree $k$ in the variable $z_n$ at $b$.
	By Proposition 1.3 we obtain that $h_C$ is polynomial in $z_n$. We see by construction and Claim 2 that $h_C$ is monic of degree $k$ in the variable $z_n$.
	\hfill$\Box_{\mathrm{Claim\,3}}$
	
	\vs{0.2cm}\noindent
	{\bf Claim 4:} $h_{C,a}$ is a Weierstrass polynomial of degree $k$ in the variable $z_n$ at every $a\in A$.
	
	\vs{0.1cm}\noindent
	{\bf Proof of Claim 4:}
	Let $a\in A$. By Claim 3, $h_a$ is polynomial of degree $k$ in the variable $z_n$. We have that $R(a)=a$.
	By Claim 2 there is an open neighbourhood $D$ of $a$ contained in $Z_C$ such that
	$h_C|_D=h^a|_D$. We are done since $h^a$ is a Weierstrass polynomial in the variable $z_n$ at $a$. 
	\hfill$\Box_{\mathrm{Claim\,4}}$
	
	\vs{0.2cm}\noindent
	So $\ma{V}:=\big\{Z_C\colon C\in \ma{C}\big\}$ and $h_V:=h_C, u_V:=u_C$ for $V=Z_C\in\ma{V}$ do the job.
\hfill$\Box$

\vs{0.5cm}
{\bf 3.4 Theorem}

\vs{0.1cm}
	Let $g\in\ma{O}^\mathrm{def}(\Omega)$. Let $k\in\big\{0,\ldots,\mathrm{ord}(f)\big\}$ and let $A$ be a definable subset of $\Omega_k$ such that $f$ is regular in the variable $z_n$ at every $a\in A$. Then there is a finite covering $\ma{V}$ of $A$ by definable open sets that are contained in $\Omega$ such that for every $V\in\ma{V}$ the following holds: 
	There are $q\in \ma{O}^{\mathrm{def}}(V)$ and  $r\in \ma{O}^{\mathrm{def}}(V^\dagger)[z_n]$ of degree at most $k-1$ in the variable $z_n$ such that
	$g=qf+r$ on $V$.
	
\vs{0.1cm}
{\bf Proof:}
	
	\vs{0.1cm}
	This follows from Proposition 2.4 with the same arguments as used in the proof of Theorem 3.3.\\
	\hfill$\Box$
	
\vs{0.5cm}
{\bf 3.5 Theorem} (Piecewise definable Weierstrass preparation)

\vs{0.1cm}	
	There is a finite covering $\ma{W}$ of $\Omega$ by definable open sets and for every $W\in\ma{W}$ there is a linear coordinate transformation $c_W$ and functions
	$h_j\in \ma{O}^{\mathrm{def}}\big(c_W^{-1}(W)^\dagger\big)[z_n]$ that is monic and $u_W\in \ma{O}^{\mathrm{def}}\big(c_W^{-1}(W)\big)$ which is a unit such that $f\circ c_W=h_Wu_W$ on $c_W^{-1}(W)$.
	
\vs{0.1cm}
{\bf Proof:}

\vs{0.1cm}
	Let $k\in\big\{0,\ldots,\mathrm{ord}(f)\big\}$. By Proposition 3.2 there is $q_k\in \IN$ and there are linear coordinate transformations $c_{k,1},\ldots,c_{k,q_k}$ such that for every $a\in \Omega_k$ there is $p\in\{1,\ldots,q_k\}$ such that
	$f\circ c_{k,p}\in\ma{O}^\mathrm{def}\big(c_{k,p}^{-1}(\Omega)\big)$ is regular of order $k$ in the variable $z_n$ at $c_{k,p}^{-1}(a)$.
	Let $B_{k,p}$ be the set of all $a\in \Omega_k$ such that $f\circ c_{k,p}$ is regular of order $k$ in the variable $z_n$ at $c_{k,p}^{-1}(a)$ and let $A_{k,p}:=c_{k,p}^{-1}(B_{k,p})$.
	By Theorem 3.3 we find a finite open covering $\ma{V}_{k,p}$ of $A_{k,p}$ by definable sets such that the following holds for every $V\in\ma{V}_{k,p}$: 
	There are
	functions 
	$h_V\in \ma{O}^{\mathrm{def}}(V^\dagger)[z_n]$ of degree $k$ in the variable $z_n$ that is monic and
	$u_V\in \ma{O}^{\mathrm{def}}(V)$ which is a unit such that
	$f\circ c_{k,p}=h_Vu_V$ on $V$.
	Take 
	$$\ma{W}:=\Big\{c_{k,p}(V)\colon V\in \ma{V}_{k,p}, k\in \big\{0,\ldots,\mathrm{ord}(f)\big\}, p\in\{1,\ldots,q_k\}\Big\}$$
	and $h_W:=h_V,u_W:=u_V$ for $W=c_{k,p}(V)\in\ma{W}$.
\hfill$\Box$

\vs{0.5cm}
{\bf 3.6 Theorem} (Piecewise definable Weierstrass division)

\vs{0.1cm}	
	Let $g\in \ma{O}^{\mathrm{def}}(V)$. There is a finite covering $\ma{W}$ by  definable open sets and for every $W\in\ma{W}$ there is a linear coordinate transformation $c_W$ and functions
	$q_W\in \ma{O}^{\mathrm{def}}\big(c_W^{-1}(W)\big)$ and  $r_W\in \ma{O}^{\mathrm{def}}\big(c_W^{-1}(W)^\dagger\big)[z_n]$ such that $g\circ c_W=q_W(f\circ c_W)+r_W$ on $c_W^{-1}(W)$.
	
\vs{0.1cm}
{\bf Proof:}

\vs{0.1cm}
	This follows from Theorem 3.4 with the same arguments as used in the proof of Theorem 3.5.
	\hfill$\Box$

\vs{0.5cm}
{\bf 3.7 Remark}
\begin{itemize}
		\item[(1)]
		We have shown  in [12] that the o-minimal structure $\IR_\an$ allows definable global complexification, meaning the following: Let $U$ be a definable open set that in $\IR^n$ and let $f\in \ma{C}_{\IR_\an}^{\omega,\mathrm{def}}(U)$. Then there is an open set $Z$ in $\IC^n$ containing $U$ that is definable and $F\in\ma{O}_{\IR_\an}^\mathrm{def}(Z)$ that extends $f$. So in the o-minimal structure $\IR_\an$ we obtain piecewise
		Weierstrass preparation and piecewise Weierstrass division for real analytic definable functions (i.e. replacing in Theorem 3.5 and Theorem 3.6
		``$\ma{O}_{\IR_\an}^\mathrm{def}$'' by ``$C^{\omega,\mathrm{def}}_{\IR_\an}$''.
		\item[(2)]
		We also obtain piecewise
		Weierstrass preparation and piecewise Weierstrass division for semialgebraic real analytic functions (= Nash functions, see Bochnak et al. [3, Chapter 8]) since we have definable global complexification in the semialgebraic setting (e.g. by [13, Section 5.2]).
		\item[(3)]
		Note that not in every o-minimal structure a definable version of the classical (i.e. local) Weierstrass theorems can be established for real analytic functions; this is for example the case with the structure $\IR_{\exp}$ (see [17, Remark 2.26] and [13, Section 3.1]).
	\end{itemize}

\newpage
\section{Definable global Nullstellensatz for principal ideals}

{\bf 4.1 Proposition}

\vs{0.1cm}
	Let $\Omega$ be a definable domain in $\IC^n$ and let $f\in \ma{O}^{\mathrm{def}}(\Omega)$ with $f\neq 0$.
	Then there is a finite covering $\ma{V}$ of $\Omega$ by definable domains such that for every $V\in\ma{V}$ the function $f|_V\in \ma{O}^{\mathrm{def}}(V)$ has a unique finite factorization into irreducible elements of $\ma{O}^{\mathrm{def}}(V)$.
	
\vs{0.1cm}
{\bf Proof:}
	
	\vs{0.1cm}
	In dimension $n=1$ the statement is obvious since then, by o-minimality, $f$ has only finitely many zeros (and one does not need the finite covering).
	The inductive step is performed using piecewise Weierstrass preparation, following the classical proof that the rings of holomorphic germs at the origin are factorial (see [11, Chapter II, Section B, Theorem 7]:)\\
	For $k\in \big\{0,\ldots,\mathrm{ord}(f)\big\}$ let $\Omega_k:=\big\{a\in\Omega\colon \mathrm{ord}_a(f)=k\big\}$. Passing to a suitable cell decomposition of $\Omega$ compatible with the $\Omega_k$'s, it is enough to show the following: Let $A$ be a definable subset of some $\Omega_k$. Then there is a finite covering $\ma{V}$ of $A$ by definable domains such that for every $V\in\ma{V}$ the function $f|_V$ has a unique factorization into irreducible elements of $\ma{O}^{\mathrm{def}}(V)$.
	By Proposition 3.2 we may assume that there is some linear coordinate transformation $c$ such that for every $a\in A$ the function $f\circ c$ is regular of order $k$ in the variable $z_n$ at $c^{-1}(a)$. Replacing $f$ by $f\circ c$ we may assume that $c$ is the identity.
	Apply Theorem 3.3 to $A$ and obtain $\ma{V}$ and for $V\in \ma{V}$ functions $h_V,u_V$ as described there. We can assume that $V$ is a domain for every $V\in\ma{V}$. 
	Fix $V\in\ma{V}$ and let $h:=h_V$. It is enough to show the claim for $h$.
	We have that $h$ is polynomial in the variable $z_n$. 
	By $V$ we denote an arbitrary element of a suitable covering of $A$ by finitely many definable domains.
	Applying the inductive hypothesis we may assume, after passing to a suitable covering of $A$ by finitely many definable domains, that the  coefficients have a unique finite factorization into irreducible elements of $\ma{O}^{\mathrm{def}}(V^\dagger)$. Following the proof of the result of Gauss that the polynomial ring over a factorial ring is factorial (see for example \L ojasiewiecz [15, A.6.2]), we obtain, after passing to a suitable covering of $A$ consisting of finitely many definable domains , that $h$ has a unique factorization into irreducible elements of $\ma{O}^{\mathrm{def}}(V^\dagger)[z_n]$.
	Applying the reasoning of the proof of [11, Chapter II, Section B, Lemma 5] and using that $h_a$ is a Weierstrass polynomial for every $a\in A$, we obtain, after passing to a suitable covering of $A$ by finitely many definable domains, that $h$ has a unique factorization into irreducible elements of $\ma{O}^{\mathrm{def}}(V)$.
\hfill$\Box$

\vs{0.5cm}
{\bf 4.2 Lemma}

\vs{0.1cm}
	Let $\Omega$ be a definable domain and let $f,g\in\ma{O}^{\mathrm{def}}(\Omega)$ with $f\neq 0$.
	The follwoing are equivalent:
	\begin{itemize}
		\item[(i)] There is a definable holomorphic function $h:\Omega\to \IC$ such that $g=hf$.
		\item[(ii)] There is a holomorphic function $h:\Omega\to \IC$ such that $g=hf$. 
	\end{itemize}
{\bf Proof:}

\vs{0.1cm}
	The implication (i)$\Rightarrow$(ii) is trivial. We show (ii)$\Rightarrow$(i):
	Let $h:\Omega\to \IC$ be the holomorphic function such that $g=fh$.
	Let $U:=\big\{z\in\Omega\colon f(z)\neq 0\}$. Then $h|_U=g|_U/f|_U$. Hence $h|_U$ is definable.
	Since $U$ is dense in $\Omega$ we have that $h$ is the continuous extension of $h|_U$ to $\Omega$. So $h$ is definable.
\hfill$\Box$

\vs{0.5cm}
\noindent
We are able to show a definable global Nullstellensatz for principal ideals.
We denote by $\ma{Z}(f)$ the zero set and by $\big(f\big)$ the principal ideal of a function $f$.

\vs{0.5cm}
{\bf 4.3 Theorem}

\vs{0.1cm}	
	Let $U$ be a definable open set in $\IC^n$ and let $f,g\in\ma{O}^{\mathrm{def}}(U)$. The following are equivalent.
	\begin{itemize}
		\item[(i)] $\ma{Z}(f)\subset\ma{Z}(g)$.
		\item[(ii)] There is some $N\in\IN$ such that $g^N\in \big(f\big)$.
	\end{itemize}
{\bf Proof:}

\vs{0.1cm}
	The implication (ii)$\Rightarrow$(i) is obvious. We show (i)$\Rightarrow$(ii), using the proof of \L ojasiewicz [15, Theorem in I.5.2]:
	By o-minimality, $U$ has finitely many connected components each of which is definable. So we may assume that $\Omega:=U$ is a domain and that $f\neq 0$.
	For $k\in \big\{0,\ldots,\mathrm{ord}(f)\big\}$ let $\Omega_k:=\big\{a\in\Omega\colon \mathrm{ord}_a(f)=k\big\}$. Passing to a suitable cell decomposition of $\Omega$ compatible with the $\Omega_k$'s it is enough to show the following: Let $A$ be a definable subset of some $\Omega_k$. Then there is a finite covering $\ma{V}$ of $A$ by definable domains such that for every $V\in\ma{V}$ there is some $M\in \IN$ such that $g|_V^M\in \big(f|_V\big)$. (Take then $N$ as the maximum of all $M=M_{A,V}$ where $A$ varies over a suitable cell decomposition.) 
	By Proposition 3.2 we may assume that there is some linear coordinate transformation $c$ such that for every $a\in A$ the function $f\circ c$ is regular of order $k$ in the variable $z_n$ at $c^{-1}(a)$. Replacing $f$ by $f\circ c$ we may assume that $c$ is the identity.
	Apply Theorem 3.3 to $A$ and obtain $\ma{V}$ and for $V\in \ma{V}$ functions $h_V,u_V$ as described there. We can assume that $V$ is a domain for every $V\in\ma{V}$. 
	Fix $V\in\ma{V}$ and let $h:=h_V$. It is enough to show the claim for $h$.
	We have that $h$ is polynomial in the variable $z_n$. 
	By $V$ we denote an arbitrary element of a suitable covering of $A$ by finitely many definable domains.
	By Proposition 4.1 we may assume, after passing to a suitable covering of $A$ consisting of finitely many definable domains, that $h$ has a unique finite factorization into irreducible elements of $\ma{O}^{\mathrm{def}}(V^\dagger)$ (use again the reasoning of the proof of [11, Chapter II, Section B, Lemma 5]).
	Hence $h$ can be written as a product $p_1^{r_1}\cdot\ldots\cdot p_s^{r_s}$ where $p_1,\ldots, p_s$ are pairwise distinct irreducible elements of $\ma{O}^{\mathrm{def}}(V^\dagger)[z_n]$ and $r_1,\ldots,r_s\in\IN$.
	Let $\hat{h}:=p_1\cdot\ldots\cdot p_s$.
	Passing to a suitable covering of $A$ by finitely many definable domains we may assume by Proposition 4.1 (and by [15, A.6.2]) that the discriminant of $\hat{h}$
	does not vanish on $V$.
	Applying Theorem 3.4 to $\hat{h}$ we can use the reasoning of the proof of [15, Theorem in I.5.2] to find some holomorphic function $\mu:V\to \IC$ such that $g|_V=\mu\hat{h}$. By Lemma 4.2 we have that $\mu\in \ma{O}^\mathrm{def}(V)$.
	Let $M:=\max\{r_1,\ldots,r_s\}.$ Then $g|_V^M\in \big(h\big)=\big(f|_V\big)$.
\hfill$\Box$

\vs{0.5cm}
\noindent
One can also use the piecewise definable Weierstrass theorems to describe definable complex analytic sets and to obtain a finite definable version of the Oka coherence result as formulated in Peterzil and Starchenko [18]. We refrain from this.

\newpage
\noi \footnotesize{\centerline{\bf References}
\begin{itemize}
\item[(1)]
G. Binyamini, D. Novikov: The Pila-Wilkie theorem for subanalytic families: a complex analytic approach.
arXiv 1605.04537.
\item[(2)]
G. Binyamini, D. Novikov: Wilkie's conjecture for restricted elementary functions.
arXiv 1605.0467.
\item[(3)]
J. Bochnak, M. Coste, M.-F. Roy: Real algebraic geometry. Ergebnisse der Mathematik und ihrer Grenzgebiete {\bf 36}, Springer, 1998.
\item[(4)]
J. Denef, L. van den Dries: $p$-adic and real subanalytic sets. {\it Ann. of Math. (2)} {\bf 128}, no. 1 (1988), 79-138.
\item[(5)]
L. van den Dries: On the elementary theory of restricted elementary functions.
{\it J. Symbolic Logic} {\bf 53}, no. 3 (1988), 796-808.
\item[(6)]
L. van den Dries: Tame Topology and O-minimal Structures. {\it London Math. Soc. Lecture Notes Series} {\bf 248}, Cambridge University Press, 1998.
\item[(7)]
L. van den Dries, C. Miller:
Extending Tamm's theorem. {\it Ann. Inst. Fourier} {\bf 44}, no. 5 (1994), 1367-1395.
\item[(8)]
L. van den Dries, P. Speissegger: The real field with convergent generalized power series. {\it Trans. Amer. Math. Soc.} {\bf 350}, no. 11 (1998),  4377-4421.
\item[(9)]
L. van den Dries, P. Speissegger: The field of reals with multisummable series and the exponential function. {\it Proc. London Math. Soc. (3)} {\bf 81}, no. 3 (2000), 513-565.
\item[(10)]
G. Fatabbi, A. Tancredi:
On the factoriality of some rings of complex Nash functions.
{\it Bull. Sci. Math.} {\bf 126}, no. 1 (2002), 61-70.
\item[(11)]
R. Gunning, H. Rossi: Analytic functions of several complex variables.
Reprint of the 1965 original. AMS Chelsea Publishing, Providence, RI, 2009.
\item[(12)]
T. Kaiser: Global complexification of real analytic globally subanalytic functions. {\it Israel Journal of Mathematics} {\bf 213} (2016),
139-174. 
\item[(13)]
T. Kaiser: $R$-analytic functions. {\it Arch. Math. Logic} {\bf 55} (2016), 605-623.
\item[(14)]
T. Kaiser, J.-P. Rolin, P. Speissegger: Transition maps at non-resonant hyperbolic singularities are o-minimal. {\it J. Reine Angew. Math.} {\it 636} (2009), 1-45.
\item[(15)]
S. \L ojasiewicz: Introduction to Complex Analytic Geometry. Birkh\"auser, 1991.
\item[(16)]
Y. Peterzil, S. Starchenko: Expansions of algebraically closed fields in o-minimal structures.
{\it Selecta Mathematica, New series} {\bf 7} (2001), 409-445.
\item[(17)]
Y. Peterzil, S. Starchenko: Expansions of algebraically closed fields II: Functions of several variables.
{\it Journal of Math. Logic} {\bf 3}, no. 1 (2003), 1-35.
\item[(18)]
Y. Peterzil, S. Starchenko: Complex analytic geometry in a nonstandard setting. Model theory with applications to algebra and analysis. Vol. 1, 117-165, {\it London Math. Soc. Lecture Note Ser.} {\bf 349}, Cambridge Univ. Press, 2008.
\item[(19)]
Y. Peterzil, S. Starchenko: Complex analytic geometry and analytic-geometric categories. {\it J. Reine Angew. Math.} {\bf 626} (2009), 39-74.
\item[(20)]
J.-C. Tougeron: Alg\`{e}bres analytiques topologiquement noeth\'{e}riennes. Th\'{e}orie de Khovanskii. {\it Ann. Inst. Fourier (Grenoble)} {\bf 41}, no. 4 (1991), 823-840.
	\end{itemize}}
	
\vs{0.5cm}
Tobias Kaiser\\
University of Passau\\
Faculty of Computer Science and Mathematics\\
tobias.kaiser@uni-passau.de\\
D-94030 Germany

\end{document}